\documentclass[a4paper,12pt]{amsart}
\usepackage{amssymb,amscd} 
\usepackage{latexsym}
\usepackage{graphicx}
\usepackage{amssymb}
\theoremstyle{plain}

\newtheorem{thm}{Theorem}[section]

\newtheorem{lem}[thm]{Lemma}
\newtheorem{cor}[thm]{Corollary}

\theoremstyle{definition}

\theoremstyle{remark}
\newtheorem{rem}{Remark}[section]
\newtheorem{exa}[rem]{Example}
\newenvironment{pf}{\medskip\noindent{Proof:}}{\hfill\qed \newline \medskip}

\title{Some families of directed strongly regular graphs
obtained from certain finite incidence structures}

\date{\today}

\newcommand{\forme}[1]{}

\newcommand{\mc}{\mathcal}

\newcommand{\vb}{\mathbf{v}}
\newcommand{\bb}{\mathbf{b}}
\newcommand{\kb}{\mathbf{k}}
\newcommand{\rb}{\mathbf{r}}

\begin{document} \pagenumbering{arabic} \setcounter{page}{1}
 \author[O.~Olmez]{Oktay Olmez}
\address{Department of Mathematics, Iowa State University,
Ames, Iowa, 50011, U. S. A.} \email[O. ~Olmez]{oolmez@iastate.edu}

\author[S. Y.~Song]{Sung Y. Song}
\address{Department of Mathematics, Iowa State University, Ames, Iowa, 50011, U. S.
A.} \email[S. Y.~Song]{sysong@iastate.edu}

\begin{abstract}

This is the second report of our work on the construction of
directed strongly regular graphs. In our previous work, we
constructed a couple of infinite families of new directed strongly
regular graphs on the sets of antiflags of partial geometries and
group divisible designs. In this paper, we use some collections of
antiflags (not the entire set of antiflags) of tactical
configurations to construct another couple of infinite families of
directed strongly regular graphs. Our construction methods are
capable of producing many, if not all, nonisomorphic directed
strongly regular graphs with same parameters.
\end{abstract}
\maketitle

\section{Introduction and preliminaries}\label{sec:intro}
The concept of directed strongly regular graphs was introduced by A.
M. Duval \cite{Du} as a generalization of the concept of strongly
regular graphs and doubly regular tournaments\footnote{A tournament
is a loopless directed graph whose adjacency matrix $A$ satisfies
$A+A^T+I=J$. A tournament is said to be doubly regular if $A$
satisfies $A^2=\lambda A+\mu(J-I-A)$ for some positive integers
$\lambda$ and $\mu$.} in 1988. The concept of strongly regular
graphs was introduced by R. C. Bose in the early 1960s, although
similar concept had been known earlier under the notion of
association schemes. The interest in strongly regular graphs has
been stimulated by the development of the theory of finite
permutation groups and the classification of finite simple groups.
It is well known that strongly regular graphs arise from many
algebraic and geometric objects including finite fields, finite
geometries, combinatorial designs and algebraic codes. The sources
for directed strongly regular graphs (with $0<t<k$) are also rich
and diverse as reported by many researchers in \cite{BH, BH2, Du,
DI, FK0, FK, HS, Jo1, Jo, KM, KP, OS}. The result of our work is to
demonstrate this claim by showing another source through several
explicit constructions.

A strongly regular graph with parameters $(v, k, \lambda, \mu)$ is
defined as an undirected regular graph $G$ with $v$ vertices
satisfying the properties that the number of common neighbors of
vertices $x$ and $y$ is $k$ if $x=y$, $\lambda$ if $x$ and $y$ are
adjacent, and $\mu$ if $x$ and $y$ are non-adjacent distinct
vertices. In terms of the adjacency matrix $A$ of a graph $G$,
identity matrix $I$ and all-ones matrix $J$, the graph $G$ is a
strongly regular graph with parameters $(v, k, \lambda, \mu)$ if and
only if (i) $JA=AJ=kJ$ and (ii) $A^2=kI+\lambda A+\mu (J-I-A)$.

A loopless directed graph $D$ with $v$ vertices is called
\textit{directed strongly regular graph} with parameters $(v, k, t,
\lambda, \mu)$ if and only if $D$ satisfies the following
conditions:
\begin{itemize}
\item[i)] Every vertex has in-degree and out-degree $k$.
\item[ii)] Every vertex $x$ has $t$ out-neighbors, all of which are
also in-neighbors of $x$.
\item[iii)] The number of directed paths of length two from a vertex
$x$ to another vertex $y$ is $\lambda$ if there is an (directed)
edge from $x$ to $y$, and is $\mu$ if there is no edge from $x$ to
$y$.\end{itemize} In terms of adjacency matrix $A=A(D)$, $D$ is a
directed strongly regular graph with parameters $(v, k, t, \lambda,
\mu)$ if and only if (i) $JA=AJ=kJ$ and (ii) $A^2=tI+\lambda A+\mu
(J-I-A)$. A strongly regular graph and a doubly regular tournament
may be viewed as a directed strongly regular graph with $t=k$ and
$t=0$, respectively. In what follows, a directed strongly regular
graph with parameters $(v, k, t, \lambda, \mu)$ will be denoted by
DSRG-$(v, k, t, \lambda, \mu)$.

In this paper, we prove existence (by explicit construction) of
directed strongly regular graphs for families of parameter sets
\begin{enumerate}
\item[(1)] $(v,k, t, \lambda,\mu) = (r(1+ab)^2b,\ r(1+ab)ab,\
ra^2b+a,\ ra^2b+a-ab-1,\ ra^2b+a)$ for any positive integers $r, a$
and $b$ such that $r\ge 2$;
\item[(2)] $(v,k, t, \lambda,\mu) = \left((1+\frac{ls}{d})s,\
ls,\ ld,\ ld-d,\ ld\right)$
and
\item[(3)] $(v,k,t,\lambda, \mu) = \left((1+\frac{ls}{d})s,\ ls+s-1,\
ld+s-1,\ ld+s-2,\ (l+1)d\right)$ for any positive integers $d, l$
and $s$ such that $d|ls$ and $1\le l< \frac{ls}{d}$.
\end{enumerate}
In this way, we confirm the existence of many graphs whose existence
was previously undetermined. Examples include the directed strongly
regular graphs with parameters $(v,k, t, \lambda, \mu)$ ($v\le 110$)
given by: $(45, 30, 22, 19, 22)$, $(54, 36, 26, 23, 26),$ $(72, 48,
34, 31, 34)$, $(75, 60, 52, 47, 52)$, $(81, 54, 38, 35, 38)$, $(90,
30, 11, 8, 11),\ (90, 60, 44, 38, 44),$ $(99, 66, 46, 43, 46)$,
$(100, 40, 18, 13, 18)$, $(108, 36, 13, 10 , 13)$, $(108, 72, 50,
47, 50),\ (108, 72, 52, 46, 52)$, and $(108, 90, 80, 74, 80)$ among
the feasible parameter sets listed in ``Parameters of directed
strongly regular graphs" by Andries Brouwer and Sylvia Hobart at
[\underline{http://homepages.cwi.nl/\~ aeb/math/dsrg/dsrg.html}].

In our construction, each directed strongly regular graph is defined
on a collection of antiflags of a tactical configuration. By
definition, a \textit{tactical configuration with parameters}
$(\vb,\bb, \kb,\rb)$ is a triple $\mc{T}=(P, \mc{B}, \mc{I})$ where
$P$ is a $\vb$-element set, $\mc{B}$ is a collection of
$\kb$-element subsets of $P$ (called `blocks') with $|\mc{B}|=\bb$,
and $\mc{I}=\{(p, B)\in P\times \mc{B}:\ p\in B\}$ such that each
element of $P$ (called a `point') belongs to exactly $\rb$ blocks.
For the notational simplicity, we will denote the tactical
configuration by the pair $\mc{T}=(P, \mc{B})$ as incidence relation
$\mc{I}$ is the natural incidence relation between the points and
blocks. When the point set and the block set are clear from the
context, we also denote a tactical configuration by
$\mc{T}-(\vb,\bb, \kb,\rb)$.

The organization of the paper is as follows. Two major construction
methods will be introduced in Section 2 and in Section 5.
Discussions on two special cases for ``Construction I" are discussed
in Section 3 and Section 4. A variation of ``Construction II" is
discussed in Section 6. We then describe thirteen nonisomorphic
graphs with parameters $(10, 4, 2, 1, 2)$ that can be obtained from
our method. The graph automorphism group of one of them acts
transitively on its vertex set; and thus, we obtain a Schurian
association scheme of class five which is a fission scheme of
Johnson scheme $J(5,2)$. We close our paper by making a few remarks
and revisiting the list of directed strongly regular graphs of small
orders found in \cite{BH} and \cite{Jo1}.


\section{Construction I: DSRG-$(r(1+ab)^2b,\ r(1+ab)ab,\ ra^2b+a,\
ra^2b+a-ab-1,\ ra^2b+a)$}\label{sec-gc}

Let $r$ and $q$ be positive integers such that $q-1=ab$ for some
positive integers $a$ and $b$. We will assume that all these
integers are greater than 1 in this section. We will consider the
case of $a=1$ and the case of $b=1$ separately in the subsequest
sections. Let $P=\{1,2,\dots, n\}$ be an $n$-element set with
$n=rq$. Let $\{G_1, G_2, \dots, G_r\}$ be a partition of $P$ into
$r$ subsets (called `groups') of size $q$. For each $j=1, 2, \dots,
r$, let $(G_j, \mc{P}_j)$ be a tactical configuration with
parameters $(\vb, \bb,\kb, \rb)=(q, q, a, a)$. That is, the block
set $\mc{P}_j$ consists of $q$ blocks such that each block is an
$a$-element subset of $G_j$, and every point of $G_j$ appears in
exactly $a$ blocks\footnote{It is easy to see that such a
configuration exists for given $q$ and $a$. For example, given a
$q$-element set $G=\{1, 2, \dots, q\}$, take $\mc{P}=\{P_1, P_2,
\dots, P_q\}$ where $P_i=\{i, i+1, i+2, \dots, i+a-1\}$ with
addition modulo $q$ to have such a configuration.}. Let the blocks
in $\mc{P}_j$ be labeled by $P_{j1}, P_{j2}, \dots, P_{jq}$. It is
clear that $P_{ig}\cap P_{jh}=\emptyset$ if $i\neq j$ since groups
are disjoint. Let $\{B_1, B_2, \dots, B_q\}$ be a family of
$ra$-element subsets of $P$ defined in such a way that
\begin{enumerate} \item[(i)] every $B_i$ contains exactly one block
from every $\mc{P}_j$, and \item[(ii)] each block in each $\mc{P}_j$
is contained in $B_i$ for exactly one $i$.
\end{enumerate}
For each $g\in G_h$, let $\{X_{g1}, X_{g2}, \dots, X_{gb}\}$ be a
partition of $G_h\setminus\{g\}$ with $|X_{gl}|=a$ for all $l\in
\{1, 2, \dots, b\}$. That is, $G_j\setminus\{g\}=X_{g1}\cup
X_{g2}\cup\cdots\cup X_{gb}$ and $|X_{gl}|=a$ for every $l$. Then we
have the following tactical configuration.
\begin{lem}
For each point $g\in G_h$, and each $l\in \{1, 2, \dots, b\}$, if we
define
\[B_{gl,j}=X_{gl}\cup (B_j\setminus P_{hj})\
\mbox{ for } j=1, 2, \dots, q,\] \[\mc{B}_g=\{B_{gl,j}:\ 1\le j\le
q,\ 1\le l\le b\},\] and
\[\mc{B}\ =\ \bigcup_{g=1}^{rq}\mc{B}_g\ =\ \{B_{gl,j}:\ 1\le g\le rq,\
1\le l\le b,\ 1\le j\le q\},\] then the pair $(P, \mc{B})$ forms a
tactical configuration with parameters
\[(\vb,\bb,\kb,\rb)=(rq,\ rq^2b,\ ra,\ rq(q-1)).\]\end{lem}
\begin{pf} From the definition, $\vb$ and $\kb$ are clear and
$\bb=|P||\mc{B}_i|=rq\cdot qb$. For $\rb$, given a point $g\in G_h$,
we have to find the size of the set $\{B'\in \mc{B}:\ g\in B'\}$. We
claim that $\rb$ is the sum of $q(q-1)$ and $(r-1)qab$. The first
summand $q(q-1)$ comes from the fact that $g$ is a member of an
$X_{il}$ for each $i\in G_h\setminus \{g\}$, and each $X_{il}$ is
contained in $q$ blocks in $\mc{B}_i$. The second summand
$(r-1)q\cdot ab$ is the number of blocks $B'$ such that $g\in B'\in
\mc{B}\setminus (\bigcup_{i\in G_h}\mc{B}_i)$ since $g$ belongs to
$B_j$ for $a$ different $j$'s (because $g$ belongs to $a$ blocks of
$(G_h, \mc{P}_h)$), and each $B_j$ is contained in $b$ blocks of
$\mc{B}_i$ for each of $(r-1)q$ points $i\in P\setminus G_h$. This
completes the proof.
\end{pf}

We now use this tactical configuration to construct a directed
strongly regular graph as follows.

\begin{thm}\label{dsrg-gtc} Let $\mc{T}$ be the above tactical configuration
$(P,\mc{B})$. Let $D=D(\mc{T})$ be the directed graph defined on the
vertex set \[V(D)=\{(g, B):\ B\in \mc{B}_g,\ g\in P\}\] with
adjacency between vertices $(g,B), (g',B')\in V(D)$ defined by $(g,
B)\rightarrow (g', B')$ if and only if $g\in B'$. Then $D$ is a
directed strongly regular graph with parameters $(v, k, t, \lambda,
\mu)$ equals to
\[(rq^2(q-1)/a, \ rq(q-1),\ r(q-1)a+a,\ q(a-1)+(r-1)(q-1)a,\
r(q-1)a+a).\]

\end{thm}

\begin{pf} It is clear that $v=|\mc{B}|=rq^2b=rq^2(q-1)/a$.
The parameter $k$ is the size of the set $\{(g',B')\in V(D):\ g\in
B'\}$ for given a vertex $(g,B)\in V(D)$, and it equals to
$\rb=rq(q-1)$. To compute $t$, let $(g,B)\in V(D)$ with $g\in G_h$
and let $B=B_{gl,j}=X_{gl}\cup (B_j\setminus P_{hj})$ for some $l$
and $j$. Then $t=|\{(g',B')\in V(D):\ g'\in B,\ g\in B'\}|$. We see
that for each $g'\in X_{gl}\subset B$, there are $q$ blocks in
$\mc{B}_{g'}$ all of which contain $g$. On the other hand, for each
of $g'\in B_j\setminus P_{hj}\subset B$, there are $ab=q-1$ blocks
in $\mc{B}_{g'}$ containing $g$. Together, we have
$t=qa+(q-1)(r-1)a$ as desired since $|X_{gl}|= a$ and $|B_j\setminus
P_{hj}|=(r-1)a$.\\

Let $(g,B)$ and $(g',B')$ be two adjacent vertices with $g\in B'$.
Suppose $g'\in G_f$ and $B'=B_{g'l, j}=X_{g'l}\cup (B_j\setminus
P_{fj})$. In order to show that $\lambda=|\{(g^*, B^*)\in V(D):\
g^*\in B',\ B^*\ni g\}|$ is constant, we consider two cases:
\begin{enumerate} \item[Case 1.] Suppose $g\in G_f$, that is, $g\in X_{g'l}$. Then (i)
for each element, say $g^*$, of $X_{g'l}\setminus \{g\}$ there are
$q$ blocks of $\mc{B}_{g^*}$ containing $g$; while (ii) for each
element $g^*\in B_j\setminus P_{fj}$, there are $ab$ blocks of
$\mc{B}_{g^*}$ containing $g$. Therefore,
$\lambda=(a-1)q+(r-1)a(q-1)$ in this case.
\item[Case 2.] If $g\notin G_f$, then $g$ must be an element of
$B_j\setminus P_{fj}$. Suppose $g\in P_{hj}\subset G_h$. Then (i)
for each choice of $g^*\in X_{g'l}$ there are $ab=q-1$ blocks
possessing $g$ (so available for $B^*$) in $\mc{B}_{g^*}$; (ii) for
each choice of $g^*\in P_{hj}\setminus \{g\}$, there are $q$ blocks
possessing $g$ in $\mc{B}_{g^*}$; and (iii) for each element $g^*$
of the remaining $(r-2)a$ elements in $B'$, there are $(q-1)$ blocks
available for $B^*$ in $\mc{B}_{g^*}$. Hence together we have
$\lambda=a(q-1)+(a-1)q+(r-2)a(q-1)$ as well.
\end{enumerate} Hence $\lambda$ has constant value
$(a-1)q+(r-1)a(q-1)$.\\

For $\mu$, let $(g,B)\nrightarrow (g',B')$, (so $g\notin B'$). Let
$g$ belong to $G_h$ for some $h$. Then by the similar counting
argument, we can verify that the number of vertices $(g^*,B^*)$ such
that $(g,B)\rightarrow (g^*,B^*)\rightarrow (g',B')$ (or
equivalently the number of choices for $g^*$ and $B^*$ such that
$g^*\in B'$ and $g\in B^*$) is $aq+(r-1)a(q-1)$ whether $g$ and $g'$
belong to the same group $G_h$ for some $h$ or not as $a$ vertices
in $B'$ can be paired with $q$ blocks while the rest can be paired
with $ab=(q-1)$ blocks. This completes the proof.
\end{pf}

As a consequence of the above theorem and the result of Duval
\cite[Theorem 7.1]{Du} on directed strongly regular graphs with
$t=\mu$, we have the following corollary.

\begin{cor} Let $r$, $q$ and $a$ be positive integers such that
$a|(q-1)$ as before. Then there exist directed strongly regular
graphs with parameters
\[(mrq^2(q-1)/a, \ mrq(q-1),\ m(rqa-ra+a),\
m\{q(a-1)+(r-1)(q-1)a\},\ m(rqa-ra+a))\] for all positive integer
$m$.
\end{cor}
\medskip

\begin{exa}\label{exa2.1} To illustrate the above construction, we consider the case when $r=2,
q=5$, and $a=b=2$. This will give us a new DSRG-$(100, 40, 18, 13,
18)$, which confirms the feasibility of the parameter set (cf.
\cite{BH}).

Let $P=\{0, 1, \dots, 9\}$, $G_1=\{1, 2, 3, 4, 5\},\ G_2=P\setminus
G_1$. $\mc{P}_1=\{12, 23, 34, 45, 15\}$ and $\mc{P}_2=\{67, 78, 89,
90, 60\}$. Then one example of tactical configuration that will
produce a DSRG-$(100, 40, 18, 13, 18)$ may be described as in the
following table. In this table entries $23$, $45$ and $2367$
represent the sets $\{2, 3\}$, $\{4, 5\}$ and $\{2, 3, 6, 7\}$
respectively\footnote{For the notational simplicity, we will remove
the brackets and commas between the elements when we list sets in a
table throughout the paper.}.

\newpage
\begin{center}{\textbf{Table \ref{exa2.1}.} The blocks $\mc{B}_i$ for each point
$i$.}\end{center} \begin{center}{
\begin{tabular}{|l|c|cccccccccc|}\hline
$i$ & $X_{i1}, X_{i2}$&$B_{i1,1}$&$B_{i1,2}$&$B_{i1,3}$&$B_{i1,4}$
&$B_{i1,5}$&$B_{i2,1}$&$B_{i2,2}$&$B_{i2,3}$&$B_{i2,4}$&$B_{i2,5}$\\
\hline
\hline 1& 23, 45&2367&2378&2389&2390&2360&4567&4578&4589&4590&4560\\
\hline 2& 13, 45&1367&1378&1389&1390&1360&4567&4578&4589&4590&4560\\
\hline 3& 12, 45&1267&1278&1289&1290&1260&4567&4578&4589&4590&4560\\
\hline 4& 12, 35&1267&1278&1289&1290&1260&3567&3578&3589&3590&3560\\
\hline 5& 12, 34&1267&1278&1289&1290&1260&3467&3478&3489&3490&3460\\
\hline 6& 78, 90&7812&7823&7834&7845&7815&9012&9023&9034&9045&9015\\
\hline 7& 89, 60&8912&8923&8934&8945&8915&6012&6023&6034&6045&6015\\
\hline 8& 79, 60&7912&7923&7934&7945&7915&6012&6023&6034&6045&6015\\
\hline 9& 67, 80&6712&6723&6734&6745&6715&8012&8023&8034&8045&8015\\
\hline 0& 67, 89&6712&6723&6734&6745&6715&8912&8923&8934&8945&8915\\
\hline
\end{tabular} }\end{center}\end{exa}

\section{Construction I $(b=1)$: DSRG-$(r(1+a)^2,\ r(1+a)a,\ ra^2+a,\
ra^2-1,\ ra^2+a)$}

Let $r$ and $q$ be positive integers greater than 1, and
$P=\{1,2,\dots, rq\}$ a set of $rq$ elements. Let $\{G_1, G_2,
\dots, G_r\}$ be a partition of $P$ into $r$ groups of size $q$. For
each $j=1, 2, \dots, r$, let $\mc{P}_j$ be the family of all
$(q-1)$-element subsets of $G_j$. Let $B_1, B_2, \dots, B_q$ be
$r(q-1)$-element subsets of $P$ defined as follows:
\begin{enumerate}\item[(1)] Select one set from each family to have
$B_1=\bigcup_{j=1}^r P_{j1}$ where $P_{j1}\in \mc{P}_j$ for $j=1, 2,
\dots, r$.
\item[(2)] For $B_2$, select one set from each $\mc{P}_j\setminus
\{P_{j1}\}$, for $j=1, 2, \dots, r$, so that $B_2=\bigcup_{j=1}^r
P_{j2}$.
\item[(3)] Continue this process to have
\[B_i=P_{1i}\cup P_{2i}\cup\cdots\cup P_{ri}\quad \mbox{where } P_{ji}\in
\mc{P}_j\setminus \{P_{j1}, P_{j2}, \dots, P_{j(i-1)}\}\] for $i=3,
4, \dots, q$.\end{enumerate} Then for each point $g\in G_h$, define
\[B_{g,j}=(G_h\setminus \{g\})\cup (B_j\setminus P_{hj})\
\mbox{ for } j=1, 2, \dots, q\] and have $\mc{B}_g=\{B_{g,1},
B_{g,2}, \dots, B_{g,q}\}$. Then with
\[\mc{B}=\bigcup_{g\in P}\mc{B}_g=\{B_{g,j}:\ 1\le g\le rq,\ 1\le
j\le q\}\] the pair $(P, \mc{B})$ becomes a tactical configuration
with parameters
\[(\vb,\bb,\kb,\rb)=(rq,\ rq^2,\ r(q-1),\ rq(q-1)).\]

\begin{thm}\label{dsrg-gtc} Let $\mc{T}$ be the above tactical configuration
$(P,\mc{B})$. Let $D=D(\mc{T})$ be the directed graph defined on the
vertex set
\[V(D)=\{(g, B):\ B\in \mc{B}_g,\ g\in P\}\] with adjacency between
vertices $(g,B)$ and $(g',B')$ defined by $(g, B)\rightarrow (g',
B')$ if and only if $g\in B'$. Then $D$ is a directed strongly
regular graph with parameters $(v, k, t, \lambda, \mu)$ equal to
\[(rq^2,\ rq(q-1),\ (q-1)(rq-r+1),\ r(q-1)^2-1,\ (q-1)(rq-r+1)).\]
\end{thm}

\begin{pf}
It is clear that $v=rq^2$, $k=q(q-1)+(r-1)q(q-1)$, and
$t=q(q-1)+(r-1)(q-1)^2$.

In order to show that $\lambda$ is constant, consider vertices
$(g,B)$ and $(g',B')$ with $(g,B)\rightarrow (g',B')$ (and so $g\in
B'$). We will consider two cases. \begin{enumerate}\item[Case 1.]
Suppose both $g$ and $g'$ belong to the same group, say $G_j$ for
some $j$. Then the number of vertices $(g^*,B^*)$ such that $B^*\ni
g$ and $g^*\in B'$ may be counted as follows. (i) Since
$G_j\setminus \{g'\}\subset B'$, with any of $q-2$ choices for $g^*$
from $G_j\setminus \{g, g'\}$, all $q$ blocks in $\mc{B}_{g^*}$
provide the legitimate pairs $(g^*, B^*)$ as every block in
$\mc{B}_{g^*}$ has $g$ in it. (ii) Since $|B'\cap (P\setminus
G_j)|=(r-1)(q-1)$, so there are $(r-1)(q-1)$ possible points
available for $g^*\in B'\cap (P\setminus G_j)$. For each point $g^*$
of these possible points, there are $q-1$ blocks possessing $g$ in
$\mc{B}_{g^*}$. Hence we must have
$\lambda=q(q-2)+(r-1)(q-1)(q-1)=r(q-1)^2 -1$.
\item[Case 2.] Suppose $g\in G_j$ for some $j$ and $g'\notin G_j$. Then
the number of ways to pick suitable $(g^*,B^*)$ may be counted as
follows: (i) With each of $q-2$ possible $g^*\in
(B'\setminus\{g\})\cap G_j$, there are $q$ blocks possessing $g$ in
$\mc{B}_{g^*}$; and thus, we can have $q(q-2)$ such vertices $(g^*,
B^*)$. (ii) With any $g^*$ of $r(q-1)-(q-1)$ points in $B'\setminus
G_j$, there are $q-1$ blocks in $\mc{B}_{g^*}$ for $B'$. Hence we
also have $\lambda =q(q-2)+(r-1)(q-1)^2$ as desired. Thus, we see
that $\lambda$ is a constant.\end{enumerate}

For $\mu$, suppose $(g,B)\nrightarrow (g',B')$, (so $g\notin B'$).
Let $g\in G_j$ for some $j$. \begin{enumerate} \item[Case 1.]
Suppose $g=g'$ and $B\neq B'$. Then vertices $(g^*,B^*)$ such that
$(g,B)\rightarrow (g^*,B^*)\rightarrow (g',B')$ may be counted as
follows. For each $g^*\in B'$, the number of blocks $B'$ in
$\mc{B}_{g^*}$ that can be paired with $g^*$ is $q$ blocks if
$g^*\in B'\cap (G_j\setminus \{g\})$, while is $q-1$ blocks if
$g^*\in B'\cap (P\setminus G_j)$. Since there are $q-1$ choices for
$g^*$ in the former and $(r-1)(q-1)$ choices for the latter, we must
have $\mu=q(q-1)+(r-1)(q-1)^2$. \item[Case 2.] If $g\neq g'$, then
$g'$ must be in $P\setminus G_j$ since neither $g$ nor $g'$ may be
in $B'$. This means $B'$ should be a block that contains all $q-1$
elements of $G_j\setminus \{g\}$. For any of $G_j\setminus \{g\}$ as
$g^*$, there are $q$ blocks that contain $g$ in $\mc{B}_{g^*}$.
(This gives us $q(q-1)$ desired vertices $(g^*,B^*)$.) For each of
$(r-1)(q-1)$ possible points in $B'\setminus G_j$, there are $(q-1)$
blocks containing $g$. Hence, we have $q(q-1)+(r-1)(q-1)^2$ for
$\mu$ in this case as well.\end{enumerate} This completes the proof.
\end{pf}

\begin{exa}\label{exa3.1} Let $r=2, q=3$, $P=\{1, 2, 3, 4, 5, 6\}$, $G_1=\{1,2,3\}$
and $G_2=\{4, 5, 6\}$. With the tactical configuration described in
Table \ref{dsrg-gtc}, we have a DSRG-$(18, 12, 10, 7, 10)$. This
graph is shown to be nonisomorphic to its orientation reversing
conjugate\footnote{By `orientation-reversing conjugate' of a graph
$D$, we mean the graph whose adjacency matrix is the transpose of
$A=A(D)$.}. By L. J{\o}rgensen \cite{Jo} we know that these are the
two nonisomorphic graphs with the parameters $(18, 12, 10, 7, 10)$
and there is no more.
\end{exa}
\begin{center}{\textbf{Table \ref{exa3.1}} $\mc{T}-(6, 18, 4, 12)$.
\qquad \textbf{Table \ref{exa3.2}} $\mc{T}-(6, 12, 3,
6)$.}\end{center}
\[\begin{array}{|l|l|}\hline
i & B_{i,j}, j=1, 2, 3\\
\hline 1 & 2356, 2346, 2345\\
\hline 2 & 1356, 1346, 1345\\
\hline 3 & 1256, 1246, 1245\\
\hline 4 & 2356, 1356, 1256\\
\hline 5 & 2346, 1346, 1246\\
\hline 6 & 2345, 1345, 1245\\ \hline
\end{array}\qquad\hskip 1in
\begin{array}{|l|l|}\hline
i & B_{i,j}, j=1, 2\\
\hline 1 & 235, 246\\
\hline 2 & 135, 146\\
\hline 3 & 415, 426\\
\hline 4 & 315, 326\\
\hline 5 & 613, 624\\
\hline 6 & 513, 524\\ \hline
\end{array}\]

\begin{exa}\label{exa3.2} Let $r=3, q=2$, $P=\{1, 2, 3, 4, 5, 6\}$, $G_1=\{1,2\}$
$G_2=\{3,4\}$ and $G_3=\{5, 6\}$. With the tactical configuration
described in Table \ref{exa3.2} above, we have a DSRG-$(12, 6, 4, 2,
4)$. It is easy to see that there are $2^6=64$ different tactical
configurations available for the given combinations of $r=3$ and
$q=2$. These 64 tactical configurations yield seven nonisomorphic
graphs. (Their adjacency matrices are given below.) It is easy to
verify that the orientation reversing conjugates, whose adjacency
matrices are the transpose of the seven adjacency matrices, are all
nonisomorphic. Therefore, our construction provides us 14 distinct
graphs with parameters $(12, 6, 4, 2, 4)$. The table showing the
description of the automorphism groups of these graphs and the size
of the isomorphism classes are followed by the adjacency matrices.

\begin{center}{\textbf{Tables 3.3} The adjacency matrices of the
graphs with parameters $(12, 6, 4, 2, 4)$ constructed in Theorem
\ref{dsrg-gtc}.}\end{center}
\[N_1=\left(\begin{array}{cccccccccccc}
0 & 0 & 1 &1 &1 &0 & 1 &0 & 1 &0 & 1 &0\\
0 & 0 & 1 &1 &1 &0 & 1 &0 & 1 &0 & 1 &0\\
1 &1 & 0 & 0 & 0 & 1 &0 & 1 &0 & 1 &0 & 1\\
1 &1 & 0 & 0 & 0 & 1 &0 & 1 &0 & 1 &0 & 1\\
1 &0 & 1 &0 & 0 & 0 & 1 & 1 &1 &0 & 1 & 0\\
1 &0 & 1 &0 & 0 & 0 & 1 & 1 &1 &0 & 1 & 0\\
0 &1 & 0 & 1 &1 &1 &0 & 0 & 0 & 1 &0 & 1\\
0 & 1 &0 & 1 &1 &1 &0 & 0 & 0 & 1 &0 & 1\\
1 &0 & 1 &0 & 1 &0 & 1 &0 & 0 & 0 & 1 &1\\
1 &0 & 1 &0 & 1 &0 & 1 &0 & 0 & 0 & 1 &1\\
0 & 1 &0 & 1 &0 & 1 &0 & 1 &1 &1 &0 & 0\\
0 &1 &0 & 1 &0 & 1 &0 & 1 &1 &1 &0 & 0\\ \end{array}\right) \qquad
N_2=\left(\begin{array}{cccccccccccc}
0 & 0 & 1 &1 &1 &0 & 1 &0 & 1 &0 & 1 &0\\
0 & 0 & 1 &1 &1 &0 & 1 &0 & 1 &0 & 1 &0\\
1 &1 &0 & 0 & 0 & 1 &0 & 1 &0 & 1 &0 & 1\\
1 &1 &0 & 0 & 0 & 1 &0 & 1 &0 & 1 &0 & 1\\
1 &0 & 1 &0 & 0 & 0 & 1 &1 &1 &0 & 0 & 1\\
1 &0 & 1 &0 & 0 & 0 & 1 &1 &1 &0 & 0 & 1\\
0 & 1 &0 & 1 &1 &1 &0 & 0 & 0 & 1 &1 &0\\
0 & 1 &0 & 1 &1 &1 &0 & 0 & 0 & 1 &1 &0\\
0 & 1 &1 &0 & 0 & 1 &1 &0 & 0 & 0 & 1 &1\\
0 & 1 &1 &0 & 0 & 1 &1 &0 & 0 & 0 & 1 &1\\
1 &0 & 0 & 1 &1 &0 & 0 & 1 &1 &1 &0 & 0\\
1 &0 & 0 & 1 &1 &0 & 0 & 1 &1 &1 &0 & 0\\ \end{array}\right)\]

\[N_3=\left(\begin{array}{cccccccccccc}
0 & 0 & 1 &1 &1 &0 & 1 &0 & 1 & 0 & 1 &0\\
0 & 0 & 1 &1 &1 &0 & 1 &0 & 1 & 0 & 1 &0\\
1 &1 &0 & 0 & 0 & 1 &0 & 1 &0 & 1 &0 & 1\\
1 &1 &0 & 0 & 0 & 1 &0 & 1 &0 & 1 &0 & 1\\
1 &0 & 1 &0 & 0 & 0 & 1 &1 &0 & 1 &1 &0\\
1 &0 & 1 &0 & 0 & 0 &1 &1 &0 & 1 &1 &0\\
0 & 1 &0 & 1 &1 &1 &0 & 0 & 1 &0 & 0 & 1\\
0 & 1 &0 & 1 &1 &1 &0 & 0 & 1 &0 & 0 & 1\\
1 &0 & 1 &0 & 0 & 1 &1 &0 & 0 & 0 & 1 &1\\
1 &0 & 1 &0 & 0 & 1 &1 &0 & 0 & 0 & 1 &1\\
0 & 1 &0 & 1 &1 &0 & 0 & 1 &1 &1 &0 & 0\\
0 & 1 &0 & 1 &1 &0 & 0 & 1 &1 &1 &0 & 0\\
\end{array}\right)\qquad
N_4=\left(\begin{array}{cccccccccccc}
0 & 0 & 1 &1 &1 &0 & 1 &0 & 1 &0 &1 &0\\
0 & 0 & 1 &1 &1 &0 & 1 & 0 &1 &0 &1 &0\\
1 &1 &0 & 0 & 0 &1 &0 & 1 &0 & 1 &0 & 1\\
1 &1 &0 & 0 & 0 &1 &0 & 1 &0 & 1 &0 & 1\\
1 &0 & 1 &0 & 0 & 0 &1 &1 &1 &0 & 1 &0\\
1 &0 & 1 &0 & 0 & 0 &1 &1 &1 &0 & 1 &0\\
0 & 1 &0 & 1 &1 &1 &0 &0 & 0 & 1 &0 & 1\\
0 & 1 &0 & 1 &1 &1 &0 &0 & 0 & 1 &0 & 1\\
1 &0 & 0 & 1 &0 &1 &1 &0 & 0 & 0 & 1 &1\\
1 &0 & 0 & 1 &0 &1 &1 &0 & 0 & 0 & 1 &1\\
0 & 1 &1 &0 & 1 &0 & 0 & 1 &1 &1 &0 & 0\\
0 & 1 &1 &0 & 1 &0 & 0 & 1 &1 &1 &0 & 0\\
\end{array}\right) \]

\[N_5=\left(\begin{array}{cccccccccccc}
0 & 0 & 1 &1 &1 &0 & 1 & 0 & 1 &0 & 1 &0\\
0 & 0 & 1 &1 &1 &0 & 1 & 0 & 1 &0 & 1 &0\\
1 &1 &0 & 0 & 0 & 1 &0 &1 &0 & 1 &0 & 1\\
1 &1 &0 & 0 & 0 & 1 &0 &1 &0 & 1 &0 & 1\\
1 &0 & 1 &0 & 0 & 0 & 1 &1 &0 & 1 &0 & 1\\
1 &0 & 1 &0 & 0 & 0 & 1 &1 &0 & 1 &0 & 1\\
0 & 1 &0 & 1 &1 &1 &0 & 0& 1 &0 & 1 &0\\
0 & 1 &0 & 1 &1 &1 &0 & 0& 1 &0 & 1 &0\\
1 &0 & 1 &0 & 1 &0 & 1 &0& 0 & 0 & 1 &1\\
1 &0 & 1 &0 & 1 &0 & 1 &0 &0 & 0 & 1 &1\\
0 & 1 &0 & 1 &0 & 1 &0 & 1 &1 &1 &0 & 0\\
0 & 1 &0 & 1 &0 & 1 &0 & 1 &1 &1 &0 & 0\\
\end{array}\right) \qquad
N_6=\left(\begin{array}{cccccccccccc}
0 & 0 & 1 &1 &1 &0 & 1 &0 & 1 &0 & 1 &0\\
0 & 0 & 1 &1 &1 &0 & 1 &0 & 1 &0 & 1 &0\\
1 &1 &0 & 0 & 0 & 1 &0 &1 &0 & 1 &0 & 1\\
1 &1 &0 & 0 & 0 & 1 &0 &1 &0 & 1 &0 & 1\\
1 &0 & 1 &0 & 0 & 0 & 1 &1 &0 & 1 &1 &0\\
1 &0 & 1 &0 & 0 & 0 & 1 &1 &0 & 1 &1 &0\\
0 & 1 &0 & 1 &1 &1 &0 & 0 & 1 &0 & 0 & 1\\
0 & 1 &0 & 1 &1 &1 &0 & 0 & 1 &0 & 0 & 1\\
0 & 1 &1 &0 & 0 & 1 &1 &0 & 0 & 0 & 1 &1\\
0 & 1 &1 &0 & 0 & 1 &1 &0 & 0 & 0 & 1 &1\\
1 &0 & 0 & 1 &1 &0 & 0 & 1 &1 &1 &0 & 0\\
1 &0 & 0 & 1 &1 &0 & 0 & 1 &1 &1 &0 & 0\\
\end{array}\right)\]

\[N_7=\left(\begin{array}{cccccccccccc}
0 & 0 & 1 &1 &1 &0 & 1 &0 & 1 &0 & 1 &0\\
0 & 0 & 1 &1 &1 &0 & 1 & 0& 1 &0 & 1 &0\\
1 &1 &0 & 0 & 0 & 1 &0 & 1&0 & 1 &0 & 1\\
1 &1 &0 & 0 & 0 & 1 &0 &1 &0 & 1 &0 & 1\\
1 &0 & 1 &0 & 0 & 0 & 1 &1 &0 & 1 &0 & 1\\
1 &0 & 1 &0 & 0 & 0 & 1 &1 &0 & 1 &0 & 1\\
0 & 1 &0 & 1 &1 & 1 &0 &0 & 1 &0 & 1 &0\\
0 & 1 &0 & 1 &1 & 1 &0 &0 & 1 &0 & 1 &0\\
0 & 1 &1 &0 & 1 &0 & 1 &0 & 0 & 0 & 1 &1\\
0 & 1 &1 &0 & 1 &0 & 1 &0 & 0 & 0 & 1 &1\\
1 &0 & 0 & 1 &0 & 1 &0 & 1 &1 &1 &0 & 0\\
1 &0 & 0 & 1 &0 & 1 &0 & 1 &1 &1 &0 & 0\\
\end{array}\right)\]

\[\begin{array}{|c|c|c|c|c|c|c|c|}\hline
\mbox{Graph} & N_1 & N_2 & N_3 & N_4 & N_5 & N_6 & N_7\\ \hline
\mbox{Automorphism Group}& D_{12}& D_8 & C_2\times C_2& C_2\times C_2 & D_{12}& S_4 & C_2\\
\hline \mbox{Size of Isomorphism Class}& 4 & 6 &12 & 12& 4& 2& 24 \\
\hline\end{array}\]
\end{exa}

\medskip
\begin{rem} In Example \ref{exa3.2}, we produced fourteen directed strongly
regular graphs with parameters $(12, 6, 4, 2, 4)$. However,
J{\o}rgensen has shown that there exist exactly twenty nonisomorphic
graphs with parameters $(12, 5, 3, 2, 2)$, which are the
complementary graphs of directed strongly regular graphs with
parameters $(12, 6, 4, 2,4)$. Therefore, there are six graphs that
are not obtained from the above construction.
\end{rem}

\begin{rem} Due to the above construction, the following feasible parameter
sets listed on the table in \cite{BH} with $v\le 110$, are realized.
\begin{center}{\textbf{Table 3.4} The new DSRGs (with $v\le 110$)
 constructed by Theorem \ref{dsrg-gtc}.}\end{center}
\begin{center}{\begin{tabular}{|c|c|ccccc|c|}\hline
 $q$ & $r$ & $v$ & $k$ & $t$ & $\lambda$ & $\mu$ & Remarks\\
\hline  3& 5& 45& 30& 22& 19& 22&\\
\hline  3& 6& 54& 36& 26& 23& 26&\\  
\hline  3& 8& 72& 48& 34& 31& 34&\\
\hline  3& 9& 81& 54& 38& 35& 38&\\
\hline  3& 11& 99& 66& 46& 43& 46&\\
\hline  3& 12& 108& 72& 50& 47& 50&\\
\hline  5 & 3& 75& 60& 52& 47& 52& \\
\hline  6& 3 & 108& 90& 80& 74& 80&\\
\hline 3& 5& 90& 60& 44& 38& 44 & m=2 (\cite[7.1]{Du})\\
\hline  3& 6& 108& 72& 52& 46& 52& m=2 (\cite[7.1]{Du})\\
 \hline
\end{tabular}}\end{center}
\end{rem}

%
\section{Construction I $(a=1)$: DSRG-$(r(1+b)^2b,\ r(1+b)b,\ rb+1,\
rb-b,\ rb+1)$}

In this section we introduce a construction method that may be
considered as a particular case of the construction method
introduced in Section \ref{sec-gc}. Although the graphs constructed
in this section may be obtained from Section \ref{sec-gc}, we
describe the construction in a different way to demonstrate the
connection to the graphs constructed in \cite{OS}. This construction
produces much more graphs than the method reported in \cite{OS}
including an infinite family of new graphs which were previously
unknown. For example, DSRG-$(90, 30, 11, 8, 11)$ for $r=5$ and
$q=3$, and DSRG-$(108, 36, 13, 10, 13)$ for $r=6$ and $q=3$ are the
new graphs among the unknown graphs listed in \cite{BH}.
\medskip

Let $r$ and $q$ be positive integers greater than 1 such that $r\le
q^{r-3}$, and let $P$ be a set of $rq$ elements. Let $\mc{P}=\{G_1,
G_2, \dots, G_r\}$ be a partition of $P$ into $r$ groups of size
$q$. Let
$$\mc{B}=\{B\subset P:\ |B\cap G_i|=1 \mbox{ for all } i=1,2, \dots,
r\}.$$ Then $\mc{B}$ consists of $q^r$ subsets (which will be called
`blocks') of $P$ of size $r$. For each $i\in P$, let
\[\mc{B}_i=\{B\in \mc{B}: i\in B\}.\] Then $|\mc{B}_i|=q^{r-1}$.
Let $\mc{B}_i$ be partitioned into $q^{r-2}$ parts each of which
consists of $q$ blocks such that no two blocks in the same part
share any other common point besides $i$. To be precise, let
$\mc{B}_{i,1}, \mc{B}_{i,2},\dots, \mc{B}_{i,w}$, where $w=q^{r-2}$,
denote the parts of the partition of $\mc{B}_i$, so that
$$\mc{B}_i=\bigcup_{j=1}^{w}\mc{B}_{i,j}$$ where (i) $\mc{B}_{i,j}\cap
\mc{B}_{i,h}=\emptyset$, for any distinct $j, h\in \{1, 2, \dots,
w\}$; (ii) $|\mc{B}_{i,j}|=q$, for every $j\in \{1, 2, \dots, w\}$;
and (iii) $B\cap C=\{i\}$ for any $B,C\in\mc{B}_{i,j}$ for each
$j\in \{1, 2, \dots, w\}$.

Given any injective map $\pi: \{1, 2, \dots, rq\}\rightarrow \{1, 2,
\dots, w\}$, if $g\in G_h$, let $\mc{C}_g^{\pi}$ denote the
collection of all blocks in the $\mc{B}_{i,\pi(i)}$ for all $i\in
G_h\setminus \{g\}$, and let $\mc{B}^{\pi}$ be the union of
$\mc{C}_g^{\pi}$ over all points in $P$. That is, for each given
injection $\pi$, define
\[\mc{B}^{\pi} = \bigcup_{g\in
P}\mc{C}_g^{\pi} \quad \mbox{ where }\quad
\mc{C}_g^{\pi}=\bigcup_{i\in G_h\setminus \{g\}}\mc{B}_{i,\pi(i)}\
\mbox{ for } g\in G_h.\] Then $\mc{T}^{\pi}=(P, \mc{B}^{\pi})$
becomes a tactical configuration with parameters \[(\vb, \bb, \kb,
\rb)=(rq,\ rq^2(q-1),\ r,\ rq(q-1)).\] We obtain a directed strongly
regular graph from this tactical configuration as follows.

\begin{thm}\label{dsrg-tc-pi} Let $D=D(\mc{T}^{\pi})$ be the
directed graph defined on the vertex set \[V(D)=\{(g, B):\ B\in
\mc{C}_g^{\pi},\ g\in P\}\] with adjacency defined by $(g,
B)\rightarrow (g', B')$ if and only if $g\in B'$. Then $D$ is a
directed strongly regular graph with parameters $(v, k, t, \lambda,
\mu)$ equals to
\[(rq^2(q-1),\ rq(q-1),\ rq-r+1,\ rq-r-q+1,\ rq-r+1).\]
\end{thm}
\begin{pf}
Since for each $g\in P$, $|\mc{C}_g^{\pi}|=q(q-1)$, we have
\[v=|V(D)|=\sum_{g\in P}|\mc{C}_g^{\pi}|=rq\cdot q(q-1).\]

A vertex $(g',B')$ is an out-neighbor of $(g,B)$ if $B'$ contains
$g$. There are $q$ blocks containing $g$ in $\mc{B}_{g,\pi(g)}$.
Every block $B'\in \mc{B}_{g,\pi(g)}$ can be paired with any point
besides the $r$ points of $B'$ to become a neighbor of $(g,B)$.
Hence we have $k=q\cdot (rq-1)$.

In order to count the (in and out)-neighbors of a vertex $(g,B)$, we
need to count the vertices $(g',B')$ such that $g'\in B$ and $B'\ni
g$. If $g$ belongs to $G_j$ for some $j$ and if $g'$ belongs to
$B\cap G_j$, then $g'$ can be paired with any block containing $g$
to become both (in and out)-neighbor of $(g,B)$. Any of the
remaining $r-1$ points belonging to $B$ (except $g'$) can be paired
with any of $q-1$ blocks containing $g$ (excluding the block
containing both $g$ and $g'$); and thus, we have $t=q+(r-1)(q-1)$.

Given $(g,B)\rightarrow (g',B')$, (and so $g\in B'$), the parameter
$\lambda$ counts the vertices $(g^*,B^*)\in V(D)$ such that $B^*\ni
g$ and $g^*\in B'$. There are $q-1$ choices for $B^*$ (except for
the block $B'$) and $r-1$ choices for $g^*$ in $B'$ excluding $g$;
and thus, $\lambda=(r-1)(q-1)$.

For $\mu$, let $(g,B)\nrightarrow (g',B')$, (and so $g\notin B'$).
If $g$ belongs to $G_j$ for some $j$ and $B'\cap G_j=\{g^*\}$, then
any block $B^*$ containing $g$ can be paired with $g^*$ to form a
path of length two from $(g,B)\rightarrow (g^*,B^*)\rightarrow
(g',B')$. Every other point in $B'$ can be paired with any of $q-1$
blocks containing $g$ (excluding the block containing both $g$ and
itself). Hence we have $\mu=q+(r-1)(q-1)$. This completes the proof.
\end{pf}

\section{Construction II: DSRG-$(ns,\ ls+s-1,\ ld+s-1,\ ld+s-2,\
ld+d)$ and\\ DSRG-$(ns,\ ls,\ ld,\ ld-d,\ ld)$ with $d(n-1)=ls$
}\label{S7}

Let $n, d, l$ and $s$ be positive integers such that $d(n-1)=ls$, or
equivalently, $n=1+\frac{ls}{d}$. Let $P=\{1, 2, \dots, n\}$. For
each $i\in P$, suppose there exists a tactical configuration
$\mc{P}_i=(P\setminus \{i\}, \mc{B}_i)$ with parameters $(\vb, \bb,
\kb, \rb)=(n-1, s, l, d)$.\footnote{There exists such a tactical
configuration subject to the conditions: $d|ls,\ 1\le d< s,\ 1\le
l\le n-2$.} We define the tactical configuration $\mc{T}=(P,
\mc{B})$ with $\mc{B}=\bigcup_{i=1}^n \mc{B}_i$ by collecting the
blocks of all configurations $\mc{P}_1, \mc{P}_2, \dots, \mc{P}_n$.
Then $\mc{T}=(P, \mc{B})$ has parameters $(\vb, \bb, \kb, \rb)=(n,
ns, l, ls)$.  Using this configuration, we now construct two
directed strongly regular graphs on the set
\[V=\{(g,B):\ B\in \mc{B}_g,\ g\in P\}.\]

\begin{thm}\label{dsrg-tcg} Let $\mc{T}=(P, \mc{B})$ be the
above tactical configuration $\mc{T}-(n, ns, l, ls)$ where
$n=1+\frac{ls}{d}$. Let $D_1=D_1(\mc{T})$ be the directed graph with
its vertex set
\[V=\{(g,B):\ B\in \mc{B}_g,\ g\in P\}\] and adjacency defined by
\[(g, B)\rightarrow (g', B')\ \mbox{ if and only if }\ g\in B'.\]
Then $D_1$ is a directed strongly regular graph with the parameters
\[(v,k,t,\lambda, \mu)=\left(ns,\ ls,\ ld,\ (l-1)d,\ ld\right).\]\end{thm}
\begin{pf} It is clear that $v=\sum_{g\in P}|\mc{B}_g|=ns$. A vertex
$(g',B')$ is to be an out-neighbor of $(g,B)$, $g'$ can be any point
different from $g$, and $B'$ can be any member of $\mc{B}_{g'}$
containing $g$. Since there are $d$ blocks in $\mc{B}_{g'}$
containing $g$, $k=(n-1)d$. A vertex $(g',B')$ is to be a (in and
out)-neighbor of $(g,B)$, $g'$ should be one of $l=|B|$ points while
$B'$ must be any one of $d$ blocks containing $g$ and belonging to
$\mc{B}_{g'}$. Hence $t=ld$.

Given $(g,B)\rightarrow (g',B')$, (and so $g\in B'$), the number of
vertices $(g^*,B^*)\in V(D)$ such that $g^*\in B'$, $B^*\in
\mc{B}_{g^*}$ and $B^*\ni g$, is $(l-1)d$ since there are $l-1$
choices for $g^*$ in $B'\setminus \{g\}$ and for any $g^*$, there
are $d$ blocks in $\mc{B}_{g^*}$ that contain $g$. Thus,
$\lambda=(l-1)d$.

For $\mu$, let $(g,B)\nrightarrow (g',B')$, (and so $g\notin B'$).
For any point $g^*$ in $B'$, there are $d$ blocks in $\mc{B}_{g^*}$
that contain $g$. Hence we have $\mu=ld$. This completes the proof.
\end{pf}

\begin{thm}\label{dsrg-tcg2} Let $\mc{T}=(P, \mc{B})$ be the same tactical
configuration as in the above theorem. Let $D_2=D_2(\mc{T})$ be the
directed graph with its vertex set
\[V=\{(g,B):\ B\in \mc{B}_g,\ g\in P\}\] and adjacency defined by
\[(g, B)\rightarrow (g', B')\ \mbox{ if and only if either } g\in B'
\mbox{ or } g=g' \mbox{ and } B\neq B'.\] Then $D_2$ is a directed
strongly regular graph with the parameters
\[(v,k,t,\lambda, \mu)=(ns,\ ls+s-1,\ ld+s-1,\ ld+s-2,\ (l+1)d).\]
\end{thm}
\begin{pf} It can be proved by the routine counting argument. \end{pf}

\begin{cor}\label{dsrg:tcg3} Let $\mc{T}=(P, \mc{B})$ be the tactical
configuration, and let $D_1=D_1(\mc{T})$ and $D_2=D_2(\mc{T})$ as in
the above theorems. In the constructions for $D_1$ and $D_2$, if we
take the multi-set consisting of $m$ copies of the vertex set $V$ as
its vertex set, we can obtain the directed strongly regular graphs
with parameters
\[(v,k,t,\lambda, \mu)=(m(ns),\ mls,\ mld,\ m(l-1)d,\ mld)\]
and \[(m(ns),\ m(ls+s)-1,\ m(ld+s)-1,\ m(ld+s)-2,\ m(l+1)d),\]
respectively.
\end{cor}
\begin{pf} Similar to the proof of (cf. \cite[7.1]{Du}). \end{pf}

\section{Construction II $(d=1)$: DSRG-$(ls^2+s,\ ls+s-1,\ l+s-1,\ l+s-2,\ l+1)$ and
\\ DSRG-$(ls^2+s,\ ls,\ l,\ l-1,\ l)$}\label{S4}

In this section, as a concrete realization of the construction
method discussed in the previous section, we consider the particular
case with $d=1$. Let $l$ and $s$ be positive integers. Consider the
$(ls+1)$-element set $P=\{1, 2, \dots, ls+1\}$. For each $i\in P$,
let $\mc{B}_i=\{B_{i1}, B_{i2}, \dots, B_{is}\}$ be a partition of
$P\setminus \{i\}$ into $s$ parts (blocks) of equal size $l$. Let
$$\mc{B}=\bigcup_{i=1}^{ls+1}\mc{B}_i=\{B_{ig}:\ 1\le g\le s,\ 1\le i\le ls+1\}.$$
Then the pair $(P, \mc{B})$ forms a tactical configuration
$\mc{T}-(ls+1, s(ls+1), l, ls)$. We construct directed strongly
regular graphs on the set
\[V=\{(i, B):\ B\in \mc{B}_i,\ i\in P\}\]
in two ways.

\begin{thm}\label{dsrg-tc} Let $(P, \mc{B})$ be $\mc{T}-(ls+1, s(ls+1), l, ls)$.
Let $D_1=D_1(\mc{T})$ be the directed graph with its vertex set
\[V=\{(i, B_{ig})\in P\times \mc{B}: 1\le i\le ls+1, 1\le g\le
s\}\] and adjacency defined by
\[(i, B_{ig})\rightarrow (j,B_{jh})\ \mbox{ if and only if }\ i\in B_{jh}.\]
Then $D_1$ is a directed strongly regular graph with the parameters
\[(v,k,t,\lambda, \mu)=(ls^2+s,\ ls,\ l,\ l-1,\ l).\]\end{thm}
\begin{pf} Straightforward. \end{pf}
\begin{thm}\label{dsrg-tc2} Let $(P, \mc{B})$ be $\mc{T}-(ls+1, ls^2+s, l, ls)$.
Let $D_2=D_2(\mc{T})$ be the directed graph with its vertex set
\[V=\{(i, B_{ig})\in P\times \mc{B}: 1\le i\le ls+1, 1\le g\le
s\}\] and adjacency defined by
\[(i, B_{ig})\rightarrow (j, B_{jh})\ \mbox{ if and only if either }
\ i\in B_{jh}\mbox{ or } i=j \mbox{ and } B_{ig}\neq B_{jh}.\] Then
$D_2$ is a directed strongly regular graph with the parameters
\[(v,k,t,\lambda, \mu)=(ls^2+s,\ ls+s-1,\ l+s-1,\ l+s-2,\ l+1).\]
\end{thm}
\begin{pf} Straightforward. \end{pf}
\begin{cor}\label{dsrg:tc3} Let $(P, \mc{B})$ be the tactical
configuration $\mc{T}-(ls+1, ls^2+s, l, ls)$ as in the above. Let
$D_1=D_1(\mc{T})$ and $D_2=D_2(\mc{T})$. In the constructions for
$D_1$ and $D_2$, if we take the multi-set consisting of $m$ copies
of the vertex set $V$ as its vertex set, we can obtain the directed
strongly regular graph with parameters
\[(v,k,t,\lambda, \mu)=(m(ls^2+s),\ mls,\ ml,\ m(l-1),\ ml)\]
and \[(m(ls^2+s),\ m(ls+s)-1,\ m(l+s)-1,\ m(l+s)-2,\ m(l+1)),\]
respectively.
\end{cor}
\begin{pf} Omitted.
\end{pf}

In the above constructions in Theorems \ref{dsrg-tc} and
\ref{dsrg-tc2}, different tactical configurations coming from
different partitions of $P$ may produce nonisomorphic graphs with
the same parameters as before. For example, for $l=s=2$, we obtain
13 different directed strongly regular graphs with the same
parameter set $(v,k,t,\lambda,\mu)=(10, 4, 2, 1, 2)$. To illustrate
the above claim and to show the connections to other combinatorial
structures, we will describe them in detail in the remainder of the
current section.

\subsection{Isomorphism classes of DSRG-$(10, 4, 2, 1, 2)$}

When $l=s=2$, the number of ways to form tactical configurations
with parameters $(\vb, \bb, \kb, \rb)=(5, 10, 2, 4)$ is $243$. Let
$\mathbf{F}$ be the set of these tactical configurations. Each
tactical configuration $\mc{T}=(P, \mc{B})\in \mathbf{F}$ gives rise
to a directed strongly regular graph $D(\mc{T})$ with its vertex set
$V(\mc{T})=\{(i, B_{ij}): i\in P,\ B_{ij}\in \mc{B}\}$ by Theorem
\ref{dsrg-tc}. Consider the action of $S_5$ on $\mathbf{F}$ under
the rule that $\mc{T}_1^{\sigma}=\mc{T}_2$ if and only if
$V(\mc{T}_1)^{\sigma}=V(\mc{T}_2)$ where
\[V(\mc{T})^{\sigma}=\{(i^{\sigma}, (B_{ij})^{\sigma}): i\in P,\
B_{ij}\in \mc{B}\}\] with natural action on $B_{ij}$; i.e.,
$(B_{ij})^\sigma =\{x^\sigma, y^\sigma\}$ if $B_{ij}=\{x, y\}$.
Under this action $\mathbf{F}$ is partitioned into seven orbits. The
tactical configurations belong to the same orbit produce isomorphic
directed strongly regular graphs. Let $\mc{T}_1$, $\mc{T}_2$,
$\dots$, $\mc{T}_7$ denote the representatives of the orbits. The
block sets of these representatives are given in Table 6.1.

\begin{center}{\textbf{Table 6.1} The block sets of the representatives of
seven orbits.}\end{center}
\[\begin{array}{|l|c|c|c|c|c|c|c|}\hline
i & \mc{B}(\mc{T}_1) & \mc{B}(\mc{T}_2) & \mc{B}(\mc{T}_3)
&\mc{B}(\mc{T}_4) &
\mc{B}(\mc{T}_5) & \mc{B}(\mc{T}_6) & \mc{B}(\mc{T}_7)\\
\hline
1 & 23,45& 23, 45& 23, 45&23, 45 &23,45& 23, 45& 23, 45\\
2 & 13,45& 13, 45& 13, 45&14, 35 &13,45&13,45&13,45 \\
3 & 12,45&14,25&12,45&15,24 &14,25 &14,25&14,25\\
4 & 12,35&12,35&12,35&13,25 &12,35 &12,35&13,25\\
5 & 12, 34&12,34&13,24&12,34 &14,23 &13,24&14,23\\ \hline
\end{array}\]\\

The Table 6.2 shows the group structure of each stabilizer of
$\mc{T}_i$, $i=1, 2, \dots, 7$ and its generators. The last row of
the table indicates the size of the orbit represented by the
corresponding tactical configuration.\\

\begin{center}{\textbf{Table 6.2} Stabilizers and the size of orbits
for the action of $S_5$ on $\mathbf{F}$.}\end{center}
\[\begin{array}{|c|c|c|c|c|c|c|}\hline
 \mc{T}_1 & \mc{T}_2 & \mc{T}_3 & \mc{T}_4 & \mc{T}_5 & \mc{T}_6 & \mc{T}_7\\ \hline
 D_8 & C_2\times C_2 & C_2 & C_5\rtimes C_4 & C_2 & C_2 & D_{10}\\
(1524), (15)(24) & (12)(45), (15)(24) &
(23)(45) & (15234), (1345)& (15)(23) & (15)(34) & (12435), (12)(45)\\
\hline
 15 & 30 & 60 & 6 & 60 & 60 & 12\\
\hline
\end{array}\]\\

Let $D(\mc{T}_i)$, $i=1, 2, \dots, 7$ be the directed strongly
regular graphs with parameters $(10, 4, 2, 1, 2)$ obtained from the
seven orbit representatives given in Table 6.1 by Theorem
\ref{dsrg-tc}. Then it is shown that the orientation-reversing
conjugates of $D(\mc{T}_i)$ for $i=1, 2, \dots, 6$ are nonisomorphic
to any of the seven. The graph $D(\mc{T}_7)$\footnote{This graph was
constructed in \cite[Ex. 4.2]{KM}.} is isomorphic to its
orientation-reversing conjugate. Therefore, together with their
conjugates, our construction produces thirteen directed strongly
regular graphs for the given parameter set. However, by J{\o}rgensen
\cite{Jo} it is known that there are sixteen graphs for the given
parameter set\footnote{Jorgensen kindly provided us the adjacency
matrices of all sixteen graphs.}.\\

\begin{center}{The adjacency matrices for seven graphs, $D(\mc{T}_1)$,
$D(\mc{T}_2), \dots$, $D(\mc{T}_7)$ are as follows.\\
(The rows of the matrices are indexed by the vertices of
corresponding graphs.)}\end{center}

\medskip
\begin{center}{\textbf{ $D(\mc{T}_1)$\hskip 2in  $D(\mc{T}_2)$}}\end{center}
\[\begin{array}{r}(1,23)\\
(2,13)\\ (3,12)\\ (3,45)\\ (4,35)\\
(5,34)\\ (1,45)\\ (4,12)\\ (2,45)\\ (5,12) \end{array}
\left(\begin{array}{rrr|rrr|rrrr}
0 & 1 & 1 & 0 & 0 & 0 & 0 & 1 & 0 & 1 \\
1 & 0 & 1 & 0 & 0 & 0 & 0 & 1 & 0 & 1 \\
1 & 1 & 0 & 0 & 1 & 1 & 0 & 0 & 0 & 0 \\ \hline
1 & 1 & 0 & 0 & 1 & 1 & 0 & 0 & 0 & 0 \\
0 & 0 & 0 & 1 & 0 & 1 & 1 & 0 & 1 & 0 \\
0 & 0 & 0 & 1 & 1 & 0 & 1 & 0 & 1 & 0 \\ \hline
0 & 1 & 1 & 0 & 0 & 0 & 0 & 1 & 0 & 1 \\
0 & 0 & 0 & 1 & 0 & 1 & 1 & 0 & 1 & 0 \\
1 & 0 & 1 & 0 & 0 & 0 & 0 & 1 & 0 & 1 \\
0 & 0 & 0 & 1 & 1 & 0 & 1 & 0 & 1 & 0
\end{array}\right)\qquad\quad
\begin{array}{r}(1,45)\\
(4,12)\\ (2,45)\\ (5,12)\\ (1,23)\\ (2,13)\\
(3,25)\\ (5,34)\\ (4,35)\\ (3,14) \end{array}
 \left(\begin{array}{rrrr|rrrrrr}
0 & 1 & 0 & 1 & 0 & 1 & 0 & 0 & 0 & 1 \\
1 & 0 & 1 & 0 & 0 & 0 & 0 & 1 & 0 & 1 \\
0 & 1 & 0 & 1 & 1 & 0 & 1 & 0 & 0 & 0 \\
1 & 0 & 1 & 0 & 0 & 0 & 1 & 0 & 1 & 0 \\ \hline
0 & 1 & 0 & 1 & 0 & 1 & 0 & 0 & 0 & 1 \\
0 & 1 & 0 & 1 & 1 & 0 & 1 & 0 & 0 & 0 \\
0 & 0 & 0 & 0 & 1 & 1 & 0 & 1 & 1 & 0 \\
1 & 0 & 1 & 0 & 0 & 0 & 1 & 0 & 1 & 0 \\
1 & 0 & 1 & 0 & 0 & 0 & 0 & 1 & 0 & 1 \\
0 & 0 & 0 & 0 & 1 & 1 & 0 & 1 & 1 & 0
\end{array}\right)\]
\newpage

\begin{center}{\textbf{ $D(\mc{T}_3)$ \hskip 2in
 $D(\mc{T}_4)$}}\end{center}
\[\begin{array}{r}(1,23)\\ (2,13)\\ (3,12)\\
(1,45)\\ (4,12)\\ (2,45)\\ (5,24)\\ (4,35)\\
(3,45)\\ (5,13)\end{array}
 \left(\begin{array}{rrr|rrrrrrr}
0 & 1 & 1 & 0 & 1 & 0 & 0 & 0 & 0 & 1 \\
1 & 0 & 1 & 0 & 1 & 0 & 1 & 0 & 0 & 0 \\
1 & 1 & 0 & 0 & 0 & 0 & 0 & 1 & 0 & 1 \\ \hline
0 & 1 & 1 & 0 & 1 & 0 & 0 & 0 & 0 & 1 \\
0 & 0 & 0 & 1 & 0 & 1 & 1 & 0 & 1 & 0 \\
1 & 0 & 1 & 0 & 1 & 0 & 1 & 0 & 0 & 0 \\
0 & 0 & 0 & 1 & 0 & 1 & 0 & 1 & 1 & 0 \\
0 & 0 & 0 & 1 & 0 & 1 & 1 & 0 & 1 & 0 \\
1 & 1 & 0 & 0 & 0 & 0 & 0 & 1 & 0 & 1 \\
0 & 0 & 0 & 1 & 0 & 1 & 0 & 1 & 1 & 0
\end{array}\right)\qquad\quad
\begin{array}{r} (1,23)\\ (2,14)\\ (4,25)\\
(5,34)\\ (3,15)\\ (1,45)\\ (4,13)\\ (3,24)\\
(2,35)\\ (5,12) \end{array}
 \left(\begin{array}{rrrrr|rrrrr}
0 & 1 & 0 & 0 & 1 & 0 & 1 & 0 & 0 & 1 \\
1 & 0 & 1 & 0 & 0 & 0 & 0 & 1 & 0 & 1 \\
0 & 1 & 0 & 1 & 0 & 1 & 0 & 1 & 0 & 0 \\
0 & 0 & 1 & 0 & 1 & 1 & 0 & 0 & 1 & 0 \\
1 & 0 & 0 & 1 & 0 & 0 & 1 & 0 & 1 & 0 \\ \hline
0 & 1 & 0 & 0 & 1 & 0 & 1 & 0 & 0 & 1 \\
0 & 1 & 0 & 1 & 0 & 1 & 0 & 1 & 0 & 0 \\
1 & 0 & 0 & 1 & 0 & 0 & 1 & 0 & 1 & 0 \\
1 & 0 & 1 & 0 & 0 & 0 & 0 & 1 & 0 & 1 \\
0 & 0 & 1 & 0 & 1 & 1 & 0 & 0 & 1 & 0
\end{array}\right)\]
\medskip

\begin{center}{\textbf{ $D(\mc{T}_5)$ \hskip 2in
 $D(\mc{T}_6)$}}\end{center}
\[\begin{array}{r} (1,23)\\ (2,13)\\ (3,25)\\
(5,23)\\ (2,45)\\ (4,12)\\ (1,45)\\ (5,14)\\
(4,35)\\ (3,14)\end{array}
 \left(\begin{array}{rrrrrrrrrr}
0 & 1 & 0 & 0 & 0 & 1 & 0 & 1 & 0 & 1 \\
1 & 0 & 1 & 1 & 0 & 1 & 0 & 0 & 0 & 0 \\
1 & 1 & 0 & 1 & 0 & 0 & 0 & 0 & 1 & 0 \\
0 & 0 & 1 & 0 & 1 & 0 & 1 & 0 & 1 & 0 \\
1 & 0 & 1 & 1 & 0 & 1 & 0 & 0 & 0 & 0 \\
0 & 0 & 0 & 0 & 1 & 0 & 1 & 1 & 0 & 1 \\
0 & 1 & 0 & 0 & 0 & 1 & 0 & 1 & 0 & 1 \\
0 & 0 & 1 & 0 & 1 & 0 & 1 & 0 & 1 & 0 \\
0 & 0 & 0 & 0 & 1 & 0 & 1 & 1 & 0 & 1 \\
1 & 1 & 0 & 1 & 0 & 0 & 0 & 0 & 1 & 0
\end{array}\right)\qquad\quad
\begin{array}{r} (1,23)\\ (2,13)\\ (3,25)\\
(5,13)\\ (1,45)\\ (4,12)\\ (2,45)\\ (5,24)\\
(4,35)\\ (3,14)\end{array}
 \left(\begin{array}{rrrrrrrrrr}
0 & 1 & 0 & 1 & 0 & 1 & 0 & 0 & 0 & 1 \\
1 & 0 & 1 & 0 & 0 & 1 & 0 & 1 & 0 & 0 \\
1 & 1 & 0 & 1 & 0 & 0 & 0 & 0 & 1 & 0 \\
0 & 0 & 1 & 0 & 1 & 0 & 1 & 0 & 1 & 0 \\
0 & 1 & 0 & 1 & 0 & 1 & 0 & 0 & 0 & 1 \\
0 & 0 & 0 & 0 & 1 & 0 & 1 & 1 & 0 & 1 \\
1 & 0 & 1 & 0 & 0 & 1 & 0 & 1 & 0 & 0 \\
0 & 0 & 1 & 0 & 1 & 0 & 1 & 0 & 1 & 0 \\
0 & 0 & 0 & 0 & 1 & 0 & 1 & 1 & 0 & 1 \\
1 & 1 & 0 & 1 & 0 & 0 & 0 & 0 & 1 & 0
\end{array}\right)\]

\medskip
\begin{center}{\textbf{ $D(\mc{T}_7)$}}\end{center}
\[\begin{array}{r} (1,23)\\ (2,13)\\ (3,25)\\
(5,23)\\ (2,45)\\ (4,25)\\ (5,14)\\ (1,45)\\
(4,13)\\ (3,14)\end{array}
 \left(\begin{array}{rrrrrrrrrr}
0 & 1 & 0 & 0 & 0 & 0 & 1 & 0 & 1 & 1 \\
1 & 0 & 1 & 1 & 0 & 1 & 0 & 0 & 0 & 0 \\
1 & 1 & 0 & 1 & 0 & 0 & 0 & 0 & 1 & 0 \\
0 & 0 & 1 & 0 & 1 & 1 & 0 & 1 & 0 & 0 \\
1 & 0 & 1 & 1 & 0 & 1 & 0 & 0 & 0 & 0 \\
0 & 0 & 0 & 0 & 1 & 0 & 1 & 1 & 0 & 1 \\
0 & 0 & 1 & 0 & 1 & 1 & 0 & 1 & 0 & 0 \\
0 & 1 & 0 & 0 & 0 & 0 & 1 & 0 & 1 & 1 \\
0 & 0 & 0 & 0 & 1 & 0 & 1 & 1 & 0 & 1 \\
1 & 1 & 0 & 1 & 0 & 0 & 0 & 0 & 1 & 0
\end{array}\right)\qquad 
\]\\

\bigskip

\begin{rem} We now describe the three graphs that are not produced by our
construction. Their adjacency matrices are given by $J_8$ and $J_9$
below, and the transpose of $J_8$ gives for the third. The graph of
$J_9$ is self-transpose and has the trivial automorphism group. The
automorphism groups for the graphs of $J_8$ and its transpose are
isomorphic to $C_2$.\end{rem}

\newpage

\begin{center}{$J_8$ \hskip 3in  $J_9$}\end{center}
\[\left(\begin{array}{rrrrrrrrrr}
0&1&1&1&1&0&0&0&0&0\\
1&0&1&1&0&0&0&1&0&0\\
1&1&0&0&1&0&0&0&0&1\\
0&0&0&0&0&1&1&1&1&0\\
0&0&0&0&0&1&1&0&1&1\\
1&1&0&0&1&0&0&1&0&0\\
1&0&1&1&0&0&0&0&0&1\\
0&0&0&0&0&1&1&0&1&1\\
0&1&1&1&1&0&0&0&0&0\\
0&0&0&0&0&1&1&1&1&0\end{array}\right)\qquad\qquad\qquad
\left(\begin{array}{rrrrrrrrrr}
0&1&1&1&1&0&0&0&0&0\\
1&0&1&1&0&0&0&1&0&0\\
1&1&0&0&1&0&0&1&0&0\\
0&0&0&0&0&1&1&0&1&1\\
0&0&0&0&0&1&1&0&1&1\\
1&1&0&0&1&0&0&1&0&0\\
1&0&1&1&0&0&0&1&0&0\\
0&0&0&0&0&1&1&0&1&1\\
0&1&1&1&0&0&0&0&0&1\\
0&0&0&0&1&1&1&0&1&0 \end{array}\right)\]
\medskip

\begin{rem} For each graph $D(\mc{T}_i)$, $i=1, 2, \dots, 7$, we can
see that the full automorphism group of $D(\mc{T}_i)$ is determined
by the stabilizer of $\mc{T}_i$ under the action of the symmetric
group $S_5$ on $\mathbf{F}$. Hence from the knowledge of the orbits
and/or the stabilizers of the permutation action of $S_5$ on
$\mathbf{F}$, we can obtain the information on the number of
distinct graphs produced by our construction. This makes us to find
exact numbers or at least some lower bounds of the number of the
isomorphism classes among the graphs constructed by our method. For
instance, as we have seen it in Table 6.2 we have the following six
different tactical configurations all of which produce graph
$D(\mc{T}_4)$.
\end{rem}

%


\begin{center}{\textbf{Table 6.3} The block sets of six tactical
configurations that are isomorphic to $\mc{T}_4$. \\ (Top row
indicates the isomorphism $\sigma\in S_5$ to the first tactical
configuration.)}\end{center}

\begin{center}{\begin{tabular}{|c|c|c|c|c|c|}\hline  (1) &(23), (45) &  (14),
(35)& (15), (24)& (13), (25)&(12), (34)\\ \hline
 23 45& 23 45& 25 34& 25 34& 24 35& 24 35\\
 14 35& 15 34& 14 35& 13 45& 15 34& 13 45\\
 15 24& 14 25& 12 45& 15 24& 12 45& 14 25\\
 13 25& 12 35& 15 23& 12 35& 13 25& 15 23\\
 12 34& 13 24& 13 24& 14 23& 14 23& 12 34\\
\hline
\end{tabular}}\end{center}

\medskip
%
%
%

\noindent Therefore, the graph $D(\mc{T}_4)$ is isomorphic to the
graphs obtained from the following vertex sets.

\[\begin{array}{|c|c|c|c|c|c|}\hline
V=V_1 &V_2& V_3& V_4& V_5& V_6\\
(1) &(23),(45) &  (14), (35)& (15), (24)& (13), (25)&(12),(34)
\\ \hline

{(1,23)}, {(1,45)}&{(1,23)}, {(1,45)}& (1,25), (1,34)&(1,25),(1,34)&
(1,24),(1,35)&(1,24), (1,35)\\

(2,14), (2,35)&{(2,15)}, {(2,34)}& (2,14),(2,35)&(2,13), (2,45)&
{(2,15)},{(2,34)}&(2,13), (2,45)\\

(3,15), (3,24)&{(3,14)}, {(3,25)}&(3,12),
(3,45)&(3,15),(3,24)&(3,12), (3,45)&{(3,14)}, {(3,25)}\\

(4,13), (4,25)&{(4,12)}, {(4,35)}&(4,15),
(4,23)&{(4,12)},{(4,35)}&(4,13), (4,25)&(4,15), (4,23)\\

(5,12), (5,34)&{(5,13)}, {(5,24)}& {(5,13)},
{(5,24)}&(5,14), (5,23)&(5,14), (5,23)&(5,12), (5,34)\\
\hline
\end{array}\]\\

\subsection{Association schemes and an SRG arising
from a DSRG-$(10, 4, 2, 1, 2)$.} Let $A$ be the adjacency matrix of
$D(\mc{T}_4)$\footnote{This graph was constructed in \cite[Sec.
5]{Du} and \cite{KM}.}, and let $\bar{A}$ be the matrix given by
\[\bar{A}_{ij}=\left \{ \begin{array}{rl}
1 & \mbox{ either }A_{ij}=1 \mbox{ or } A_{ji}=1,\\
0& \mbox{ otherwise}.\\ \end{array}\right . \]

\[\bar{A} = \ \begin{array}{r} (1,23)\\ (2,14)\\ (4,25)\\
(5,34)\\ (3,15)\\ (1,45)\\ (4,13)\\ (3,24)\\
(2,35)\\ (5,12) \end{array}
 \left(\begin{array}{rrrrr|rrrrr}
0 & 1 & 0 & 0 & 1 & 0 & 1 & 1 & 1 & 1 \\
1 & 0 & 1 & 0 & 0 & 1 & 1 & 1 & 0 & 1 \\
0 & 1 & 0 & 1 & 0 & 1 & 0 & 1 & 1 & 1 \\
0 & 0 & 1 & 0 & 1 & 1 & 1 & 1 & 1 & 0 \\
1 & 0 & 0 & 1 & 0 & 1 & 1 & 0 & 1 & 1 \\ \hline
0 & 1 & 1 & 1 & 1 & 0 & 1 & 0 & 0 & 1 \\
1 & 1 & 0 & 1 & 1 & 1 & 0 & 1 & 0 & 0 \\
1 & 1 & 1 & 1 & 0 & 0 & 1 & 0 & 1 & 0 \\
1 & 0 & 1 & 1 & 1 & 0 & 0 & 1 & 0 & 1 \\
1 & 1 & 1 & 0 & 1 & 1 & 0 & 0 & 1 & 0
\end{array}\right)\]\\

Let $G$ be the graph whose adjacency matrix is $\bar{A}$. Then $G$
is the strongly regular graph with parameters $(v, k, \lambda,
\mu)=(10, 6, 3, 4)$, which is known as Johnson graph $J(5,2)$. We
note that $J(5,2)$ is also obtained from the J{\o}rgensen's graph of
$J_9$ by `symmetrizing' the matrix $J_9$. In fact, this is the only
strongly regular graph that can be obtained from any of the directed
strongly regular graphs with parameters $(10, 4, 2, 1, 2)$ through
the symmetrization process.

Among the directed strongly regular graphs with parameters $(10, 4,
2, 1, 2)$, $D(\mc{T}_4)$ has the largest automorphism group. It is
the only one that has vertex transitive automorphism group. The
automorphism group $H=\mbox{Aut}(D(\mc{T}_4))$ is isomorphic to the
group $C_5\rtimes C_4$ of order $20$. From the transitive
permutation group $H$ on the vertex set of $D(\mc{T}_4)$, we obtain
a 5-class association scheme. Let $\mc{X}(H, V(\mc{T}_4))$ denote
this association scheme. Then its association relation table is
given by the matrix on the left below.
\medskip
\begin{center}{\textbf{Tables 6.4} Relation matrices of
 $\mc{X}(H, V(\mc{T}_4))$ and its 2-class symmetric fusion scheme.}\end{center}
\[
\left(\begin{array}{ccccc|ccccc}
0 &  3 &  2 &  2 & 3 & 5 &  1 &  4 &  4 &  1\\
3 &  0 &  3 &  2 & 2 & 4 &  4 &  1 &  5 &  1\\
2 &  3 &  0 &  3 & 2 & 1 &  5 &  1 &  4 &  4\\
2 &  2 &  3 &  0 & 3 & 1 &  4 &  4 &  1 &  5\\
3 &  2 &  2 &  3 & 0 & 4 &  1 &  5 &  1 &  4\\
\hline
5 &  1 &  4 &  4 & 1 & 0 &  3 &  2 &  2 &  3\\
4 &  1 &  5 &  1 & 4 & 3 &  0 &  3 &  2 &  2\\
1 &  4 &  4 &  1 & 5 & 2 &  3 &  0 &  3 &  2\\
1 &  5 &  1 &  4 & 4 & 2 &  2 &  3 &  0 &  3\\
4 &  4 &  1 &  5 & 1 & 3 &  2 &  2 &  3 &  0\\
\end{array}\right)\qquad\qquad
 \left(\begin{array}{rrrrr|rrrrr}
0 & 1 & 2 & 2 & 1 & 2 & 1 & 1 & 1 & 1 \\
1 & 0 & 1 & 2 & 2 & 1 & 1 & 1 & 2 & 1 \\
2 & 1 & 0 & 1 & 2 & 1 & 2 & 1 & 1 & 1 \\
2 & 2 & 1 & 0 & 1 & 1 & 1 & 1 & 1 & 2 \\
1 & 2 & 2 & 1 & 0 & 1 & 1 & 2 & 1 & 1 \\ \hline
2 & 1 & 1 & 1 & 1 & 0 & 1 & 2 & 2 & 1 \\
1 & 1 & 2 & 1 & 1 & 1 & 0 & 1 & 2 & 2 \\
1 & 1 & 1 & 1 & 2 & 2 & 1 & 0 & 1 & 2 \\
1 & 2 & 1 & 1 & 1 & 2 & 2 & 1 & 0 & 1 \\
1 & 1 & 1 & 2 & 1 & 1 & 2 & 2 & 1 & 0
\end{array}\right)\]\\

It is observed that $\mc{X}(H, V(\mc{T}_4))$ is isomorphic to the
5-class non-commutative association scheme labeled as $\mathcal{Y}$
on \cite[p.255]{SS}. This scheme has two fusion schemes of class 3;
they are $K_2\times K_5$, the direct product of two trivial schemes
of order 2 and 5, and $C_5\wr K_2$, the wreath product of the scheme
coming from pentagon and the trivial scheme of order 2. The scheme
$\mc{X}(H, V(\mc{T}_4))$ also has three symmetric fusion schemes of
class 2, $K_2\wr K_5, K_5\wr K_2$ and the Johnson scheme $J(5,2)$.
The Johnson scheme $J(5,2)$ is obtained from $\mc{X}(H,
V(\mc{T}_4))$ by fusing the relations $R_1, R_3$, and $R_4$
together, and fusing $R_2$ and $R_5$ together as easily observed
from the above relation tables.

The edge set of $D(\mc{T}_4)$ coincides with $R_1\cup R_3$. The
orientation-reversing conjugate of $D(\mc{T}_4)$ is the graph with
edge set $R_1\cup R_4$. The edge set of Johnson graph $J(5, 2)$ is
$R_1\cup R_3\cup R_4$, while its complement, the Petersen graph has
edge set $R_2\cup R_5$. As we have mentioned earlier, although both
graphs $D(\mc{T}_4)$ and $J_9$ give rise to Johnson graph $J(5,2)$
via the symmetrization process, $D(\mc{T}_4)$ is the one which
yields $\mc{X}(H, V(\mc{T}_4))$.

\section{Concluding remarks}

\begin{rem} We have seen that each of the construction methods discussed
in this paper is capable of producing many different directed
strongly regular graphs with same parameters for some parameter
sets. The number of distinct graphs depends on the number of
isomorphism classes of the underlying tactical configurations. For
example, by the method discussed in Theorem \ref{dsrg-tc}, we
enumerate at least 1985 nonisomorphic DSRG-$(14, 6, 3, 2, 3)$ and
$217194772$ nonisomorphic graphs with parameters $(18, 8, 4, 3, 4)$.
However, our construction can not generate all the graphs with a
given parameter set in general. As it is known by J{\o}rgensen,
there must be 16495 graphs with parameters $(14, 6, 3, 2, 3)$. This
makes it clear that we need more work to find all graphs.

We also note that the determination of the graph automorphisms is
deduced to the investigation of the isomorphism classes of the
underlying tactical configurations as different tactical
configurations may produce isomorphic graphs. In the examples
discussed in the last section, we have seen that the number of
nonisomorphic graphs is determined by the orbit structure of the
permutation group $S_{ls+1}$ on the set of all tactical
configurations $\mc{T}-(ls+1, s(ls+1), l, ls)$ for given $l$ and
$s$. Although it is involved as the order of a graph gets large, it
is routine to calculate the automorphism groups.
\end{rem}

\begin{rem} All our constructions are based on tactical
configurations which arise in many structures. For example,
interesting particular cases of the construction methods in Theorem
\ref{dsrg-tc} and Theorem \ref{dsrg-tc2} occur when we consider a
$2-(\vb,\kb, 1)$ design, especially, a $2-(n^2+n+1, n+1, 1)$ design,
the symmetric design coming from a projective plane of order $n$.
Let $P$ be the point set of this projective plane. For a point $p\in
P$, let $L_{p0}, L_{p1}, \dots, L_{pn}$ denote the $n+1$ lines
passing through $p$. Since $p$ is the unique common intersecting
point for any two of these lines, if we set $B_{pi}=L_{pi}-\{p\}$
for $i=0, 1, \dots, n$, then with $\mc{B}=\{B_{pi}: p\in P,\ i\in
\{0, 1, \dots, n\}\}$, the pair $(P,\mc{B})$ forms a tactical
configuration with parameters $(\vb, \bb, \kb, \rb)=(n^2+n+1,
(n+1)(n^2+n+1), n, n(n+1))$. For each prime power $n$, using this
tactical configuration, we can obtain a directed strongly regular
graphs with parameters \[(v,k,t,\lambda, \mu)=((n+1)(n^2+n+1),\
n(n+1),\ n,\ n-1,\ n)\] and \[(v,k,t,\lambda, \mu)=((n+1)(n^2+n+1),\
n(n+2),\ 2n,\ 2n-1,\ n+1).\] For these directed strongly regular
graphs, we have complete information on their automorphism groups
from the knowledge of automorphism groups of projective planes. It
will be interesting to know, in what extent, any properties of the
graphs involved shed a light in the study of the geometry and vice
versa.\end{rem}

\begin{rem} Finally, we close our paper by revisiting the table of small
directed strongly regular graphs $(v\le 20)$ provided by Brouwer and
Hobart in \cite{BH} and J{\o}rgensen in \cite{Jo1} to recall the
current status of their existence, enumeration results, known
construction methods and our constructions. For undefined symbols in
the table we refer the readers to \cite{BH}. (See also the tables
provided in \cite{Du} and \cite{KM} for other characteristics of
some of these graphs.)\end{rem}

\newpage

{\small{
\begin{center}{\textbf{Table} The list of small directed strongly regular graphs
 ($v\le 20$) revisited\footnote{For some parameter sets, only one of a complementary pair
of graphs is listed in this table.}.}\end{center}

\begin{center}{\begin{tabular}{|ccccc|l|rl|rr|}\hline
  &   &   &         &     & Known construction & New &construction & Known \# & \\
 $v$ &$k$&$t$&$\lambda$&$\mu$& recorded in \cite{BH} & with  & $d, m, s, r, q, l$& \cite{Jo1} &\\ \hline
6 &2 & 1 & 0& 1& T1, T5, T8, T12 & Theorem \ref{dsrg-tc} & $ m=1, s=2, l=1$ &1 & \\
  &3 & 2 & 1& 2& T8, T9 & Theorem \ref{dsrg-tc2} & $m=1, s=2, l=1$& 1& \\ \hline
8 &3 & 2 & 1& 1&T4, T6, T7 & & &1 & \\
  &4 & 3 & 1& 3& T17 & Theorem \ref{dsrg-tc-pi} & $m=1, q=2, r=2$& 1& \\ \hline
10&4 & 2 & 1& 2& T3, T5, T12, M1 & Theorem \ref{dsrg-tc} & $m=1, s=2, l=2$&16 & \\
  &5 & 3 & 2& 3& & Theorem \ref{dsrg-tc2} & $m=1, s=2, l=2$&16 & \\ \hline
12&3 & 1 & 0& 1& T1, T8, T12 &Theorem \ref{dsrg-tc} & $m=1, s=3, l=1$& 1& \\
  &8 & 6 & 5& 6& T8 & Theorem \ref{dsrg-tcg2} & $d=2, s=3, l=2$&1 & \\
  \hline
12&4 & 2 & 0& 2& T8, T10, T12& Theorem  \ref{dsrg-tc} & $m=2, s=2, l=1$&1 & \\
  &7 & 5 & 4& 4& T8, T11 & Theorem  \ref{dsrg-tc2} & $m=2, s=2,
  l=1$& 1& \\ \hline
12&5 & 3 & 2& 2& T4, T6, T8, T9, T11& Theorem  \ref{dsrg-tc2} & $m=1, s=3, l=1$&20 & \\
  &6 & 4 & 2& 4& T8, T10, T12 & Theorem  \ref{dsrg-tc-pi} & $m=1, q=2, r=3$&20 & \\
  \hline
14&5 & 4 & 1& 2& DNE by \cite{KM}&   & & & \\
\hline
14&6 & 3 & 2& 3& T5, T12, M6&  Theorem  \ref{dsrg-tc} & $m=1, s=2, l=3$&16495 & \\
  &7 & 4 & 3& 4&  & Theorem  \ref{dsrg-tc2} & $m=1, s=2, l=3$&16495 & \\ \hline
15&4 & 2 & 1& 1& T2, T4&  & & 5& \\
\hline
15&  5& 2&   1 &  2& M5&  && 1292& \\
\hline
16&  6&   3&   1&   3 & DNE by \cite{FK0}&  & & & \\
\hline
16&  7&   4&   3&   3& T4, T6, T15& && & \\
  &  8&  5& 3& 5&  & Theorem \ref{dsrg-tc-pi} & $m=2, q=2, r=2$& & \\
  \hline
16&  7&   5&   4&   2&T11 &  && 1& \\
  &  8&   6&   2&   6& T10&  Theorem \ref{dsrg-tc-pi} & $m=1, q=2, r=4$&1 & \\
  \hline
18&  4&   3&   0&   1& M3 & && 1& \\
\hline
18&  5&   3 &  2&   1& T7& && 2& \\
\hline
18&  6&   3 &  0&   3 &T8, T10, T12, T17& Theorem \ref{dsrg-tc} & $m=3, s=2, l=1$&1 & \\
  &  11&  8 &  7&   6 & T8, T11& Theorem \ref{dsrg-tc2} & $m=3, s=2,
  l=1$& 1& \\ \hline
18&  7 &  5 &  2&   3& T16, M4& && & \\
\hline
18&  8 &  4 &  3&   4&T3, T5, T12& Theorem \ref{dsrg-tc} & $m=1, s=2, l=4$& & \\
  &  9 &  5 &  4&   5& &Theorem \ref{dsrg-tc2} & $m=1, s=2, l=4$& & \\
  \hline
18&  8 &  5 &  4&   3& T11&  && & \\
  &  9 &  6 &  3&   6& T10 & & & & \\ \hline
20&  4 &  1 &  0&   1& T1, T8, T12& Theorem \ref{dsrg-tc} & $m=1, s=4, l=1$& 1& \\
  &  15&  12&  11&  12& T8&  & & 1& \\ \hline
20&  7 &  4 &  3 &  2 &T8, T9&  Theorem \ref{dsrg-tc2} & $m=1, s=4, l=1$& & \\
  &  12&  9 &  6 &  9&T8, T12 && & & \\ \hline
20&  8 &  4 &  2 &  4& T10, T12& Theorem  \ref{dsrg-tc} & $m=2, s=2, l=2$& & \\
  &  11&  7 &  6 &  6& T11 & Theorem  \ref{dsrg-tc2} & $m=2, s=2,
  l=2$& & \\ \hline
20&  9 &  5 &  4 &  4& T4, T6, T11& && & \\
  &  10&  6 &  4 &  6 & & Theorem \ref{dsrg-tc-pi} & $m=1, q=2, r=5$& & \\
  \hline
\end{tabular}}\end{center}}}

\end{document}